\documentclass{amsart}
\usepackage{amscd}
\usepackage{thmdefs}

\theoremstyle{definition}
\theoremstyle{remark}
\numberwithin{equation}{section}

\begin{document}
\title[]{On continuous but differentiable nowhere solutions in the moment
problem}
\author{Miguel Tierz}
\address{Institut d'Estudis Espacials de Catalunya (IEEC/CSIC). Campus UAB, Fac
Ciencies, Torre C5-Parell-2a planta\\
E-08193 Bellaterra (Barcelona) Spain}
\email{tierz@ieec.fcr.es}
\date{}


\maketitle


Since the times of Chebyshev and, especially, after Stieltjes seminal memoir 
\cite{Sti}, the moment problem has always been a subject of remarkable
interest. Its notorious influence and appearance in many areas is well-known 
\cite{Akh,ShoTam,Sim}.

Let us briefly present the fundamentals of the moment problem. Let $%
I\subseteq R$ be an interval. For a positive measure $\mu $ on $I$, the $n$%
th moment is defined as $\int_{I}dxd\mu \left( x\right) $, provided the
integral exists. Consider now the set $\left\{ \gamma _{n}\right\}
_{n=0}^{\infty }$, an infinite but countable sequence of real numbers. Then,
the moment problem on $I$ consists in the solution of the following
questions:

\subsubsection{Does there exist a positive measure on I with moments $%
\left\{ \gamma _{n}\right\} _{n=0}^{\infty }$ ?}

\strut \smallskip

If so:

\subsubsection{Is this positive measure uniquely determined by the moments $%
\left\{ \gamma _{n}\right\} _{n=0}^{\infty }$ ?}

\strut \smallskip

One may also want to consider a third question, in case the previous one is
answered in the negative:

\subsubsection{How can all the positive measures on $I,$ with moments $%
\left\{ \gamma _{n}\right\} _{n=0}^{\infty },$ be described ?}

\strut \smallskip

This last question is the main focus here. The simplest way to show the
non-uniqueness feature is to recall Stieltjes famous computation\footnote{%
Stieltjes considered $k=1$.} \cite{Sti}:

\begin{equation}
\int_{0}^{\infty }x^{n}\mathrm{e}^{-k^{2}\log ^{2}x}\sin \left( 2\pi \frac{%
\log x}{\log q}\right) dx=0,\qquad n\in {\mathbb{Z},}  \label{zero}
\end{equation}
where $q=\mathrm{e}^{-1/2k^{2}}$. Thus, independently of $\lambda ,$ we have:

\begin{eqnarray}
&&\frac{1}{\sqrt{\pi }}\int_{0}^{\infty }x^{n}\mathrm{e}^{-k^{2}\log
^{2}x}\left( 1+\lambda \sin \left( 2\pi \frac{\log x}{\log q}\right) \right)
dx  \label{q} \\
&=&\frac{1}{\sqrt{\pi }}\int_{0}^{\infty }x^{n}\mathrm{e}^{-k^{2}\log
^{2}x}dx=q^{\left( n+1\right) ^{2}/2}.  \notag
\end{eqnarray}
Our main point is almost trivial but apparently absent in the literature:
instead of $\left( \ref{zero}\right) $, one can more generally write: 
\begin{equation}
\int_{0}^{\infty }x^{n}\mathrm{e}^{-k^{2}\log ^{2}x}\sin \left( 2\pi j\frac{%
\log x}{\log q}\right) dx=0,\qquad n,j\in {\mathbb{Z},}
\end{equation}
which means that one can construct the Fourier series: 
\begin{equation}
f\left( x\right) =\frac{1}{\sqrt{\pi }}\mathrm{e}^{-k^{2}\log ^{2}x}\left(
1+\lambda \sum_{n=1}^{m}a_{n}\sin \left( 2\pi b_{n}\frac{\log x}{\log q}%
\right) \right) ,  \label{logper}
\end{equation}
and the integer moments keep their (log-normal) value $\left( \ref{q}\right) 
$. Now, it is manifest that a suitable choice of $a_{n}$ and $b_{n}$ and the
case $m\rightarrow \infty $ lead to functions that are continuous but
differentiable nowhere. The sum-integral exchange in the moment computation
of $\left( \ref{logper}\right) $ (term by term integration) is, by virtue of
Lebesgue dominated convergence theorem, not a problem in the $m\rightarrow
\infty $ case, as soon as the series is convergent.

One can get a more compact and elegant description by employing a more
general characterization of the non-uniqueness, first found by Chihara in
the seventies \cite{Chihara} and more recently exploited in \cite{Christ}.
It says that the moments remain unchanged if $f\left( x\right) =\frac{1}{%
\sqrt{\pi }}k\mathrm{e}^{-k^{2}\log ^{2}x}$ is modified in the following
way: 
\begin{equation}
f\left( x\right) \rightarrow f\left( x\right) g\left( x\right) ,\text{ with }%
g\left( x\right) =g\left( qx\right) .  \label{qper}
\end{equation}

This $q$-periodic function $g\left( x\right) $ is well known in the $q$%
-calculus literature and is the most general solution of $D_{g}g\left(
x\right) =0$ where $D_{g}g\left( x\right) =\frac{g\left( x\right) -g\left(
qx\right) }{\left( q-1\right) x}$ denotes the usual $q$-derivative \cite{EE}%
. Interestingly enough, this $q$-periodic function has been studied in the
context of fractal geometry \cite{EE}. It turns out that nowhere
differentiable functions are perfectly well-behaved under the $q$-derivative
so, while in $q$-calculus a $q$-periodic function essentially plays the role
of a constant, it may easily be differentiable nowhere under the ordinary
derivative.

The property $\left( \ref{qper}\right) $ was found employing the functional
equation satisfied by the log-normal function, $f\left( xq\right) =\sqrt{q}%
xf\left( x\right) $. This functional equation actually corresponds to the $q$%
-Pearson equation associated to the Stieltjes-Wigert orthogonal polynomials.
It can be easily shown that the same argument in \cite{Christ} can be
identically applied to any orthogonal polynomials system that is an
indeterminate moment problem (all within the $q$-deformed Askey scheme \cite
{Koe}). Thus, the property $\left( \ref{qper}\right) $, also holds for all
of them and hence, the existence of differentiable nowhere solutions appears
in the same way. A detailed account of this and other related results will
appear elsewhere \cite{Tierz}.

\begin{remark}
To conclude, we would like to make a remark of historical nature. In
addition to Stieltjes seminal memoir, another relevant mathematical document
is the correspondence between Hermite and Stieltjes \cite{BaiBour}. There
are two well-known letters that are still often quoted.

One is a letter from Hermite to Stieltjes (number 351, dated 22 October
1892) which shows the admiration and respect for the results of Stieltjes.
This was a reply to previous correspondence from Stieltjes, where he
described his discoveries to him. Among these results, the one we have
discussed here. In his own words:

''L'existence de ces fonctions $\varphi \left( u\right) $ qui, sans \^{e}tre
nulles, son telles que

\begin{equation*}
\int_{0}^{\infty }u^{k}\varphi \left( u\right) du=0\qquad \left(
k=0,1,2,3,...\right) 
\end{equation*}
me para\^{i}t tr\`{e}s remarquable.''

The second one, dated 20 May (1893), is even more well-known, and is the one
where Hermite emphatically expresses his dislike for the then recently
introduced continuous but differentiable nowhere functions.
\end{remark}

Therefore, it seems interesting that the same results that impressed
Hermite, actually contained, as we have shown, mathematics that he despised
so much. Namely, continuous but differentiable nowhere functions.

\end{document}